\documentclass[a4paper,12pt]{article}
\usepackage{amsthm,amsmath}
\usepackage{amssymb}
\usepackage[all]{xypic}
\usepackage{enumerate}
\usepackage{a4wide}
\usepackage{amsthm}

\newtheorem{theorem}{Theorem}[section]

\newtheorem{corollary}[theorem]{Corollary}
\newtheorem{lemma}[theorem]{Lemma}
\newtheorem*{theorem*}{Theorem}
\newtheorem*{proposition*}{Proposition}
\newtheorem*{corollary*}{Corollary}
\newtheorem*{lemma*}{Lemma}\theoremstyle{definition}
\newtheorem{definition}[theorem]{Definition}

\newtheorem{example}[theorem]{Example}
\newtheorem{remark}[theorem]{Remark}

\newtheorem*{definition*}{Definition}

\newcommand{\coring}[1]{\mathfrak{#1}}
\newcommand{\tensor}[1]{\otimes_{#1}}
\newcommand{\rdual}[1]{#1^{\ast}}
\newcommand{\rcomatrix}[2]{#2^* \tensor{#1} #2}
\newcommand{\rend}[2]{\mathrm{End}_{-#1}({#2})}
\newcommand{\lend}[2]{\mathrm{End}_{#1 -}({#2})}
\newcommand{\rhom}[3]{\mathrm{Hom}_{-#1}({#2},{#3})}
\newcommand{\lhom}[3]{\mathrm{Hom}_{#1 -}({#2},{#3})}
\newcommand{\ldual}[1]{{ }^{\ast}#1}

\def\eps{{\varepsilon}}
\def\ot{{\otimes}}
\def\ut{{\otimes}}
\def\cC{\coring{C}}
\def\eC{\eps_\cC}
\def\DC{\Delta_\cC}
\def\cD{\coring{D}}
\def\eD{\eps_\cD}
\def\DD{\Delta_\cD}
\def\sw#1{{\sb{(#1)}}}

\begin{document}
\title{On comatrix corings and  bimodules}
\author{Tomasz Brzezi\'nski\footnote{E-mail: T.Brzezinski@swansea.ac.uk} \\ Department of Mathematics\\ University of Wales Swansea,\\ Swansea SA2 8PP, U.K.
\and Jos\'e G\'{o}mez-Torrecillas\footnote{E-mail: torrecil@ugr.es} \\ Departamento de \'Algebra\\
Universidad de Granada\\ E18071 Granada, Spain} \maketitle
\abstract{To any bimodule which is finitely generated and
projective on one side one can associate a coring, known as a
comatrix coring. A new description of comatrix corings in terms of
data reminiscent of a Morita context is given. It is also studied
how  properties of bimodules  are reflected in the associated
comatrix corings. In particular it is shown that separable
bimodules give rise to coseparable comatrix corings, while
Frobenius bimodules induce Frobenius comatrix corings.}

\date
\section{Introduction}
One of the first and most fundamental examples of corings is
provided by the canonical coring of Sweedler \cite{Swe:pre} which
can be associated to any ring extension $B\to A$. The structure of
the canonical coring detects whether such  an extension is
separable, split or Frobenius. Comodules of this coring provide
one with an equivalent description of  the descent theory for an
extension $B\to A$. In recent paper \cite{KaoGom:com} it has been
realised that Sweedler's canonical corings are special examples of
more general class of corings termed {\em comatrix corings}. A
comatrix $A$-coring can be associated to any $(B,A)$-bimodule $M$
provided $M$ is a finitely generated projective right $A$-module.
It is natural to expect that such a coring should reflect
properties of module $M$ in a way similar to the relationship
between properties of ring extensions and those of corresponding
canonical corings.

The aim of this paper it to study properties of comatrix corings
in relation to properties of bimodules. In particular we show that
the dual $(A,B)$-bimodule $\rdual{M}$ is a separable bimodule if
and only if  the corresponding comatrix coring is a cosplit
coring.  On the other hand if $M$ is a separable  (resp.\
Frobenius)  bimodule then the comatrix coring  is coseparable
(resp.\ Frobenius) coring. The converse holds provided certain
faithful flatness condition is satisfied. Since for any
$(B,A)$-bimodule $M$ one can consider a ring extension $B\to S$,
where $S$ is the right endomorphism ring of $M$, there is also
associated canonical Sweedler's coring. We study how the above
properties of a comatrix coring are reflected by the properties of
corresponding Sweedler's coring. This coring formulation of
properties of modules can shed new light on module theoretic
conjectures such as the Caenepeel-Kadison conjecture on
biseparable and Frobenius extensions \cite{CaeKad:bis} (cf.\
\cite{BrzKad:bis} for a coring formulation of the problem).

The paper is organised as follows. In Section~2 we give a new
formulation of comatrix corings in terms of algebraic data which
are very similar (semi-dual) to Morita contexts. This formulation
of comatrix corings puts them in a broader perspective of
established algebraic theories and can suggest new applications
to, for instance, K-theory. In Section~3 we study the properties
of a comatrix coring associated to a $(B,A)$-bimodule $M$  in
relation to module properties of $M$.
\medskip

\noindent{\bf  Notation and preliminaries.} Throughout the paper,
$A$ and $B$ are associative rings with 1. For modules we use the
standard module theory notation, for example a $(B,A)$-bimodule
$M$ is often denoted by ${}_BM_A$, $\rhom A--$ denotes the Abelian
group of right $A$-module maps, $\lhom B--$ denotes left
$B$-module maps etc. The left endomorphisms of $M$ are denoted by
$\lend BM$ and their ring structure is provided by opposite
composition of maps. Similarly, the right endomorphisms of $M$ are
denoted by $\rend AM$ and their ring structure is provided by
composition of maps. The dual of a right $A$-module $M$ is denoted
by $M^*$, while the dual of a left $A$-module $N$ is denoted by
$\ldual{N}$. The identity morphism of $M$ is denoted also by $M$.
Finally, for a finitely generated projective module $M_A$,
elements of a finite dual basis are denoted by $e_i\in M$ and
$\rdual{e_i}\in \rdual{M}$, $i\in I$.

An $(A,A)$-bimodule $\cC$ is called an {\em $A$-coring} provided
there exist $(A,A)$-bimodule maps $\DC:\cC\to \cC\ot_A\cC$ and
$\eC:\cC\to A$ such that
$$
(\cC\ot\DC)\circ\DC = (\DC\ot\cC)\circ\DC, \qquad
(\eC\ot\cC)\circ\DC = (\cC\ot\eC)\circ\DC = \cC.
$$
The map $\DC$ is known as a {\em coproduct} or {\em
comultiplication}, while $\eC$ is called a {\em counit}. To
indicate the action of $\DC$ we use the Sweedler sigma notation,
i.e., for all $c\in \cC$,
$$
\DC(c) = \sum c\sw1\ot c\sw 2, \qquad (\DC\ot\cC)\circ\DC(c) =
(\cC\ot\DC)\circ\DC(c) = \sum c\sw1\ot c\sw 2\ot c\sw 3,
$$
etc. A {\em morphism of $A$-corings} is an $(A,A)$-bimodule map
$\vartheta: \cC\to \cD$ such that $(\vartheta\ot\vartheta)\circ
\DC = \DD\circ\vartheta$ and $\eD\circ\vartheta = \eC$. Given any
$A$-coring $\cC$, its left dual $\ldual{\cC}$ is a ring with the
multiplication
$$
(\varphi \varphi')(c) = \sum\varphi(c\sw 1\varphi'(c\sw 2)),
\qquad \forall \varphi , \varphi'\in \ldual{\cC}, \;\; c\in\cC,
$$
and the unit $\eC$.

Given any ring extension $B\to A$, the $(A,A)$-bimodule   $\cC =
A\ot_BA$ is an $A$-coring with the coproduct $$ \DC : A\ot_BA\to
A\ot_BA\ot_AA\ot_BA\cong A\ot_BA\ot_BA, \qquad a\ot a'\mapsto a\ot
1_A\ot a',
$$
and  the counit $\eC: A\ot_BA\to A$, $a\ot a'\mapsto aa'$
\cite{Swe:pre}.  Through the natural identification $\lhom
A{A\ot_BA} A\cong \lend BA$, the left dual ring of $\cC$ is
anti-isomorphic to the endomorphism ring $\lend BA$. The coring
$\cC$ is known as a  {\em Sweedler's $A$-coring} associated to a
ring extension $B\to A$. This is the most fundamental example of a
coring, thus it is often termed a {\em canonical coring}.  For
other examples of corings and further details about their
structure and properties we refer to  \cite{Brz:str,KaoGom:com}
and to forthcoming monograph \cite{BrzWis:cor}.

\section{Comatrix corings from contexts}
Let $M$ be a  $(B,A)$-bimodule such that $M_A$ is finitely
generated and projective module. Denote by $M^* = \rhom AMA$ is
the dual $(A,B)$-bimodule and let $\{ e_i\in M, e_i^*\in \rdual{M}
\}_{i\in I}$ be a finite dual basis of $M$. Then   the
$(A,A)$-bimodule $M^*\otimes_B M$ is an  $A$-coring with the
coproduct
$$
\Delta_{M^*\otimes_B M} : M^*\otimes_B M\to M^*\otimes_B
M\otimes_AM^*\otimes_B M, \quad \varphi\otimes m\mapsto \sum_{i\in
I}\varphi\otimes e_i\otimes e_i^*\ot m,
$$
and the counit
$$
\eps_{M^*\otimes_B M} : M^*\otimes_B M\to A , \qquad
\varphi\otimes m\mapsto \varphi(m).
$$
The coring  $M^*\otimes_B M$  is known as a {\em comatrix
$A$-coring} \cite{KaoGom:com}.  Note that the definition of the
coproduct does not depend on the choice of a dual basis (cf.\
\cite[Remark~1]{KaoGom:com}). In this section we show that
comatrix corings can be understood in terms of data very
reminiscent of Morita contexts in the classical module theory.

\begin{definition}\label{def.com.context}
    Given a pair  of algebras  $A,B$, a {\em comatrix coring context}
    consists of an $(A,B)$-bimodule $N$, a
    $(B,A)$-bimodule $M$ and a pair of bimodule maps
    $$
    \sigma : N\otimes_{B}M\to A, \qquad \tau : B\to
    M\otimes_{A}N,
    $$
    such that the following diagrams
    $$
    \xymatrix{ N\ot_{B}M\ot_{A}N
    \ar[rr]^{\sigma\otimes {N}}&& A\ot_{A}N
    \ar[d]^{\cong}& & M\ot_{A}N\ot_{B}M
    \ar[d]_{ {M}\ut\sigma} && B\ot_{B}M\ar[ll]_{\tau\ut {M}}
    \ar[d]^{\cong}\\
    N\ot_{B}B\ar[u]^{{N}\ut\tau}  \ar[rr]^{\cong} && N && M\ot_{A}A \ar[rr]^{\cong} && M}
    $$
    commute. A comatrix coring context is denoted by $(A,B,{}_AN_B,{}_BM_A,\sigma,\tau)$.
\end{definition}

\begin{example}\label{ex.context.Mor}
Suppose that ${}_AN_B$ and ${}_BM_A$ together with  bimodule maps
$\sigma: N\otimes_{B}M\to A$, $\tilde{\tau} : M\otimes_{A}N\to B$
form a Morita context.   Suppose that $\tilde{\tau}$ is
surjective. By standard arguments in Morita theory (cf.\
\cite[Ch.~II.3]{Bas:KTh}) one proves that $\tilde{\tau}$ is
bijective, and let $\tau$ be the inverse of $\tilde{\tau}$. Then
$(A,B,{}_AN_B,{}_BM_A,\sigma,\tau)$ is a comatrix context.
\end{example}
\begin{proof}
This follows immediately from the definition of a Morita context.
\end{proof}

Example~\ref{ex.context.Mor} justifies the use of the term {\em
context} in Definition~\ref{def.com.context}. The following
example  explains the appearance of words {\em comatrix} and {\em
coring}.

\begin{example}\label{ex.context.comat}
Let $M$ be a  $(B,A)$-bimodule such that $M_A$ is finitely
generated and projective module with a dual basis  $\{ e_i, e_i^*
\}_{i\in I}$. Define $\tau: B\to M\otimes_A\rdual{M}$ by $b\mapsto
\sum_{i\in I} be_i\otimes e_i^* = \sum_{i\in I} e_i\otimes
e_i^*b$. Then
$$
(A,B,{}_A\rdual{M}_B,{}_BM_A,\eps_{M^*\otimes_B M},\tau)
$$
is a comatrix coring context.
\end{example}
\begin{proof}
This follows immediately from the properties of a dual basis.
\end{proof}

Example~\ref{ex.context.comat} has the following converse, which
constitutes the main result of this section.

\begin{theorem}\label{thm.com.context}
Let $(A,B,{}_AN_B,{}_BM_A,\sigma,\tau)$ be a comatrix coring
context. Let $e = \tau(1_B)\in M\otimes_A N$. Then
\begin{enumerate}[(1)]
\item $M$ is a finitely generated and projective right $A$-module and ${}_AN_B$ is isomorphic to $\rdual{M}$.
\item $\cC= N\otimes_B M$ is an $A$-coring with the coproduct
$$
\DC :\cC\to\cC\ot_A\cC, \qquad n\otimes m\mapsto n\otimes e\otimes
m,
$$
and counit $\eC = \sigma$.
\item  The coring $\cC$ is isomorphic to the comatrix coring $\rdual{M}\otimes_BM$.
\end{enumerate}
\end{theorem}
\begin{proof}
(1) Write $e = \sum_{i\in I} m_i\ot n_i$. Since $e = \tau(1_B)$,
the second of the diagrams in Definition~\ref{def.com.context}
implies that for all $m\in M$, $m = \sum_i m_i \sigma(n_i\ot m)$,
i.e., $M_A$ has a finite dual basis $\{ m_i, \sigma(n_i\ot
-)\}_{i\in I}$. Therefore $M_A$ is a finitely generated and
projective module. Furthermore, the $(A,B)$-bimodule map $\chi : N
\to \rdual{M}$, $n\mapsto \sigma(n\otimes -)$ is an
$(A,B)$-bimodule isomorphism with the inverse $\chi^{-1} :
\rdual{M}\to N$, $\varphi\mapsto \sum_i\varphi(m_i)n_i$. Indeed,
take any $n\in N$ and $m\in M$ and use the first of the diagrams
in Definition~\ref{def.com.context} to compute
$$
(\chi^{-1}\circ \chi)(n\ot m) = \sum_i \sigma(n\ot m_i)n_i\ot m =
n\ot m.
$$
Similarly take any $\varphi\in\rdual{M}$ and $n\in M$, and use the
facts that $\sigma$ is $(A,A)$-bilinear and that $\{ m_i,
\sigma(n_i\ot -)\}$ is a dual basis to compute
$$
(\chi\circ \chi^{-1})(\varphi\ot m) = \sum_i \sigma(\varphi
(m_i)n_i\ot -)\ot m = \sum_i \varphi (m_i\sigma(n_i\ot -))\ot m =
\varphi\ot m,
$$
as required.

(2) Note that $\DC$ is well-defined since the fact that $\tau$ is
a $(B,B)$-bimodule map implies that $e$ is $B$-invariant, i.e.,
$e\in (M\otimes_A N)^B = \{x\in M\otimes_AN\; |\; \forall b\in B,
\; bx=xb\}$. $\DC$ is coassociative directly from its definition.
Finally, the counit properties of $\eC = \sigma$ follow from the
commutative diagrams in Definition~\ref{def.com.context}. For
example
$$
(\cC\ot \eC)\circ\DC(n\ot m) = n \ot (M\ot \sigma)(e\ot m) = n \ot
(M\ot \sigma)(\tau(1)\ot m) = n\ot m,
$$
be the second of these diagrams. This proves that $\cC$ is an
$A$-coring.

(3)  Write $e = \sum_i m_i\ot n_i$, and let $\chi$ be the
$(A,B)$-bimodule isomorphism constructed in (1). The induced map
$\vartheta = \chi\ot M:  N\otimes_B M \to \rdual{M}\otimes_B M$,
$n\otimes m\mapsto \sigma(n\otimes -)\otimes m$ is then an
isomorphism of $(A,A)$-bimodules. Note that for all $m\in M$ and
$n\in N$, $\eps_{\rdual{M}\otimes_B M}(\vartheta(n\ot m)) =
\sigma(n\otimes m) = \eC(n\otimes m)$. Furthermore,
\begin{eqnarray*}
(\vartheta \ot \vartheta)(\DC (n\ot m)) &=& \sum_i \vartheta(n\ot m_i)\ot \vartheta( n_i\ot m) \\
&=& \sum_i \sigma(n\ot -)\ot m_i\ot \sigma(n_i\ot -) \ot m =
\Delta_{\rdual{M}\otimes_B M}(\vartheta(n\ot m)),
\end{eqnarray*}
since the definition of the coproduct in a comatrix coring does
not depend on the choice of a dual basis. Thus we conclude that
$\vartheta$ is a morphism of $A$-corings.

The inverse of $\vartheta$ is given by $\chi^{-1}\otimes_B M$ and
comes out as
$$
\vartheta^{-1} : \rdual{M}\otimes_B M\to N\ot_B M, \qquad
\varphi\ot m\mapsto \sum_i\varphi(m_i)n_i\ot m.
$$
The fact that $\sigma$ is an $(A,A)$-bimodule map and the second
of the diagrams in Definition~\ref{def.com.context} facilitate the
following calculation for all $\varphi\in \rdual{M}$ and $m\in M$
$$
\sigma(\sum_i\varphi(m_i)n_i \ot m) =
\varphi(\sum_im_i\sigma(n_i\ot m))= \varphi(m).
$$
This means that $\eC\circ \vartheta^{-1} =
\eps_{\rdual{M}\otimes_B M}$. Furthermore
\begin{eqnarray*}
(\vartheta^{-1}\ot \vartheta^{-1})(\Delta_{\rdual{M}\otimes_B
M}(\varphi\ot m)) &=&
\sum_ i \vartheta^{-1}(\varphi\ot m_i)\ot \vartheta^{-1}(\sigma(n_i\ot -)\ot m) \\
&=& \sum_{i,j,k} \varphi(m_j)n_j\ot m_i\ot \sigma(n_i\ot m_k)n_k\ot m\\
&=& \sum_{i,j} \varphi(m_j)n_j\ot m_i\ot n_i\ot m\\
&=& \DC(\vartheta^{-1}(\varphi\ot m)),
\end{eqnarray*}
where the third equality follows from the first of the diagrams in
Definition~\ref{def.com.context}. Thus $\vartheta^{-1}$ is also an
$A$-coring morphism. Consequently, $\vartheta$ is an $A$-coring
isomorphism and  we conclude that the coring $\cC$ is isomorphic
to the comatrix coring $\rdual{M}\otimes_B M$ as asserted.
\end{proof}

In view of Example~\ref{ex.context.comat},
Theorem~\ref{thm.com.context} asserts that comatrix coring
contexts provide one with an equivalent description of comatrix
corings. As an immediate consequence of
Theorem~\ref{thm.com.context} we also obtain the following
description of a left dual ring of the coring associated to a
comatrix coring context.

\begin{corollary}\label{cor.com.context}
Let $(A,B,{}_AN_B,{}_BM_A,\sigma,\tau)$ be a comatrix coring
context and let $\cC = N\ot_B M$ be the associated $A$-coring.
Then the ring $\ldual{\cC}$ is anti-isomorphic to the endomorphism
ring $\lend BM$.
\end{corollary}
\begin{proof}
Since the coring $\cC$ is isomorphic to the comatrix coring
$\rdual{M}\ot_B M$ there is a ring isomorphism $\ldual{\cC}\cong
\ldual{(\rdual{M}\ot_B M)}$. The latter ring is anti-isomorphic to
$\lend BM$ by \cite[Proposition~1]{KaoGom:com}.
\end{proof}

\section{Comatrix corings of separable and Frobenius bimodules}
This section is devoted to studies of relationship between
properties of comatrix corings and the following two notions from
the classical module theory. Let $M$ be a $(B,A)$-bimodule.
Following Sugano \cite{Sug:not} (cf.\ \cite{CaeKad:bis}), $M$ is
called a {\em separable bimodule} or $B$ is said to be {\em
$M$-separable over $A$} provided the evaluation map
$$
M\ot_A\ldual{M}\to B, \qquad m\ot\varphi\mapsto \varphi(m),
$$
is a split epimorphism of $(B,B)$-bimodules. Following
\cite{AndFul:rin} and \cite{Kad:new} a bimodule  ${}_BM_A$ is said
to be {\em Frobenius} if both ${}_BM$ and $M_A$ are finitely
generated and projective and $\ldual{M}\cong \rdual{M}$ as
$(A,B)$-bimodules.

These properties of a bimodule ${}_BM_A$ lead to corresponding
properties of the ring extension $B\to S=\rend{A}{M}$,
$b\mapsto[m\mapsto bm]$. As shown in \cite{Sug:not} (cf.\
\cite[Theorem~3.1]{Kad:sep}), if $M$ is a separable bimodule, then
$B\to S$ is a {\em split extension}, i.e., there exists a
$B$-bimodule map $\mathsf{s}:S\to B$ such that $\mathsf{s}(1_S) =
1_B$. Conversely, if ${}_BM_A$  is such that $M_A$ is finitely
generated projective, and $B\to S$ is a split extension, then
${}_BM_A$ is a separable bimodule. Furthermore, the endomorphism
ring theorem (cf. \cite[Theorem~2.5]{Kad:new}) asserts that if
${}_BM_A$ is a Frobenius bimodule, then $B\to S$ is a Frobenius
extension, i.e., $S_B$ is  finitely generated projective and
$\rdual{S}\cong S$ as $(B,S)$-bimodules.

Before we begin the discussion of the relationship of module
properties of $M$ and the properties of the corresponding comatrix
coring we make the following clarifying
\begin{remark}\label{rem.s}
Let ${}_BM_A$ be a bimodule with $M_A$ finitely generated and
projective, and consider its right endomorphism ring
$S=\rend{A}{M}$. Then there is a canonical isomorphism of
$S$-bimodules $M \tensor{A} \rdual{M} \cong S$ which sends a
simple tensor $m \tensor{} \varphi \in M \tensor{A} \rdual{M}$ to
the endomorphism $x \mapsto m\varphi(x)$. Its inverse is given by
the assignment $s \mapsto \sum _i s(e_i) \tensor{} e_i^*$,  where
 $\{ e_i, e_i^* \}_{i\in I}$ is a finite dual basis. From now on, we always identify $M
\tensor{A} \rdual{M}$ and $S$. With this identification, the
product in the ring $S$ (the composition) obeys the following
rules: given $s \in S, m \tensor{} \varphi, m' \tensor{} \varphi'
\in M \tensor{A} \rdual{M}$,
\begin{eqnarray*}
s(m \tensor{} \varphi) & = &  s(m) \tensor{} \varphi, \\
 (m
\tensor{} \varphi)s & = &  m \tensor{} \varphi s, \\
 (m \tensor{}
\varphi)(m' \tensor{} \varphi') & = & m\varphi(m') \tensor{}
\varphi' = m \tensor{} \varphi(m') \varphi'.
\end{eqnarray*}
In this case we can consider a comatrix $A$-coring $\rcomatrix
BM$. Furthermore, since there is a ring map $B\to S$, there is
also canonical Sweedler's $S$-coring $S\ot_BS$.  These are the
corings which reflect the structure of $M$, and thus they will be
of special interest in this section.
\end{remark}

An $A$-coring $\cC$ is said to be {\em cosplit} provided there
exists an $(A,A)$-bimodule section of the counit, i.e., iff $\eC$
is a split epimorphism of $A$-bimodules (cf.\
\cite[26.14]{BrzWis:cor}). The following theorem provides one with
the complete description of separability of the dual module
${}_A\rdual{M}_B$.

\begin{theorem}\label{th.cosplit}
Let ${}_BM_A$ be a bimodule such that $M_A$ is a finitely
generated projective, and let $S = \rend{A}{M}$. Then
\begin{enumerate}[(1)]
\item ${}_A\rdual{M}_B$ is a separable bimodule if and only if the
comatrix coring $\rcomatrix{B}{M}$ is a cosplit $A$--coring.
\item If the comatrix coring $\rcomatrix{B}{M}$ is a cosplit $A$-coring
then Sweedler's coring $S \tensor{B} S$ is a cosplit $S$-coring.
\end{enumerate}
\end{theorem}
\begin{proof}
$(1)$ $\rdual{M}$ is a separable bimodule if and only if the
evaluation map $\rdual{M} \tensor{B} \ldual{(\rdual{M})}
\rightarrow A$ is a split epimorphism of $A$-bimodules. Using the
natural isomorphism $\ldual{(\rdual{M})} \cong M$ the evaluation
map coincides with the counit $\eps_{\rcomatrix{B}{M}}:
\rcomatrix{B}{M} \rightarrow A$.

$(2)$ Since $\rcomatrix{B}{M}$ is a cosplit coring, there is an
$A$-bimodule map $e : A \rightarrow \rcomatrix{B}{M}$ such that
$\eps_{\rcomatrix{B}{M}} \circ e = A$. Now use the correspondence
between $S$ and $M\ot_A \rdual{M}$ discussed in Remark~\ref{rem.s}
and define $\widetilde{e} : S \rightarrow S \tensor{B} S$ as the
composite
\[
\xymatrix{S = M \tensor{A} \rdual{M} \cong M \tensor{A} A
\tensor{A} \rdual{M} \ar^-{M \ot e \ot \rdual{M}}[rr] & & M
\tensor{A} \rcomatrix{B}{M} \tensor{A} \rdual{M} = S \tensor{B}
S}.
\]
Clearly $\widetilde{e}$  is an $S$--bimodule map. We need to prove
that $\widetilde{e}$ is a splitting of the counit of the canonical
coring $S \tensor{B} S$. Recall that the counit $\eps_{S
\tensor{B} S}$   is simply the multiplication map $S \tensor{B} S
\rightarrow S$.  Write $e(1_A) = \sum_\alpha w_{\alpha}^*
\tensor{} w_{\alpha} \in \rcomatrix{B}{M}$, and note that
$\sum_\alpha w_{\alpha}^*(w_{\alpha}) = 1_A$. Identify $1_S$ with
$ \sum_i e_i\ot_A\rdual{e_i}$. Then $\widetilde{e}(1_S) =
\sum_{i,\alpha}e_i \tensor{} w_{\alpha}^* \tensor{} w_{\alpha}
\tensor{} e_i^*$ and, therefore, the multiplication map evaluated
at $\widetilde{e}(1_S)$ gives
\[
\sum_{i,\alpha}(e_i \tensor{} w_{\alpha}^*)(w_{\alpha} \tensor{}
\rdual{e_i}) = \sum_{i,\alpha}e_i \rdual{w_{\alpha}}(w_{\alpha})
\tensor{} \rdual{e_i} = \sum_i e_i \tensor{} \rdual{e_i} = 1_S.
\]
Since $\widetilde{e}$ is an $S$-bimodule map, we deduce that it
splits the counit of $S \tensor{B} S$, i.e., $S \tensor{B} S$ is a
cosplit $S$-coring.
\end{proof}

\begin{definition}\label{def.n.integral}
Given an $A$-coring $\cC$, an $A$-bimodule map $\gamma :
\coring{C} \tensor{A} \coring{C} \rightarrow A$ such that for all
$c, c' \in \coring{C}$,
$$
\sum c_{(1)} \gamma (c_{(2)} \tensor{} c') = \sum \gamma (c
\tensor{} c'_{(1)}) c'_{(2)}
$$
is called a \emph{pre-cointegral}.
\end{definition}

\begin{lemma}\label{cointegral}
If the comatrix coring $\rcomatrix{B}{M}$ has a pre-cointegral
$\gamma$, then the composite map $\widetilde{\gamma}$ given by
\[
\xymatrix{M \tensor{A} \rcomatrix{B}{M} \tensor{A}
\rcomatrix{B}{M} \tensor{A} \rdual{M} \ar^-{M \tensor{} \gamma
\tensor{} \rdual{M}}[rr] & & M \tensor{A} A \tensor{A} \rdual{M}
\cong M \tensor{A} \rdual{M}}
\]
 is a pre-cointegral
for $S \tensor{B} S$.
\end{lemma}
\begin{proof}
Note that we implicitly identify $S$ with $M\ot_A\rdual{M}$ as in
Remark~\ref{rem.s}. Obviously, $\widetilde{\gamma} : S \tensor{B}
S \tensor{B} S \rightarrow S$ is an $S$-bimodule map. Furthemore
for all $s, s', s'' \in S$ we compute
\begin{eqnarray*}
s \tensor{} \widetilde{\gamma}(1_S \tensor{} s' \tensor{} s'') & =
& \sum_{i,k} s \tensor{} \widetilde{\gamma}(e_i \tensor{}
\rdual{e_i} \tensor{}
s' \tensor{} s''(e_k) \tensor{} \rdual{e_k}) \\
 & = & \sum_{i,k} s \tensor{} e_i \gamma(\rdual{e_i} \tensor{} s' \tensor{}
 s''(e_k)) \tensor{} \rdual{e_k} \\
  & = & \sum_{i,j,k}s(e_j) \tensor{} \rdual{e_j} \tensor{} e_i \gamma(\rdual{e_i} \tensor{} s' \tensor{}
 s''(e_k)) \tensor{} \rdual{e_k} \\
  & = & \sum_{i,j,k}s(e_j) \tensor{} \gamma(\rdual{e_j} \tensor{} s' \tensor{}
  e_i)\rdual{e_i} \tensor{} s''(e_k) \tensor{} \rdual{e_k} \\
  & = & \sum_{i,j}s(e_j) \tensor{} \gamma(\rdual{e_j} \tensor{} s' \tensor{}
  e_i)\rdual{e_i} \tensor{} s'' \\
  & = & \sum_{i,j}\widetilde{\gamma}(s(e_j) \tensor{} \rdual{e_j} \tensor{}
  s' \tensor{} e_i \tensor{} \rdual{e_i}) \tensor{} s''
  =  \widetilde{\gamma}(s \tensor{} s' \tensor{} 1_S) \tensor{}
  s''.
\end{eqnarray*}
Here the identification in Remark~\ref{rem.s} has been used in
derivation of the first, third, fifth and seventh equalities. The
fourth equality follows from the fact that $\gamma$ is a
pre-cointegral. In view of the definition of a coproduct in
Sweedler's coring this proves that $\gamma$ is a pre-cointegral.
\end{proof}

A pre-cointegral $\gamma$ is called a {\em cointegral} provided
$\gamma\circ\DC = \eC$. Following  \cite{Guz:coi}, an $A$-coring
$\cC$ is called a {\em coseparable coring} provided its coproduct
$\DC$ is a split monomorphism of $(\cC,\cC)$-bicomodules.
Equivalently, by \cite[Theorem~3.5, Corollary~3.6]{Brz:str} an
$A$-coring is a coseparable coring provided it has a cointegral.
Coseparable corings turn out to correspond to separable bimodules.

\begin{theorem}\label{prop.cosep}
Let ${}_BM_A$ be a bimodule such that $M_A$ is a finitely
generated projective, and let $S = \rend{A}{M}$. Then
\begin{enumerate}[(1)]
\item If $M$ is a separable bimodule, then the comatrix coring $\rcomatrix{B}{M}$ is a
coseparable $A$--coring.
\item If a comatrix coring $\rcomatrix{B}{M}$ is a coseparable $A$--coring then $S
\tensor{B} S$ is a coseparable $S$--coring.
\end{enumerate}
\end{theorem}
\begin{proof}
$(1)$ Since ${}_BM_A$ is separable, $B \rightarrow S$ is a split
extension (cf.\ \cite[Theorem~3.1]{Kad:sep}). Let $\mathsf{s} : S
\rightarrow B$ be a $B$-bimodule splitting of the unit map. With
the identification in Remark~\ref{rem.s}, this means that
$\mathsf{s}(\sum_i e_i\ot\rdual{e_i}) = 1_B$. Define $\gamma :
\rcomatrix{B}{M} \tensor{A} \rcomatrix{B}{M} \rightarrow A$ as the
composite map
\[
\xymatrix{\rcomatrix{B}{M} \tensor{A} \rcomatrix{B}{M}
\ar^-{\rdual{M} \tensor{B} \mathsf{s} \tensor{B} M}[rr] & &
\rdual{M} \tensor{B} B \tensor{B} M \cong \rcomatrix{B}{M}
\ar^-{\eps_{\rcomatrix{B}{M}}}[rr] & &A. }
\]
Clearly, $\gamma$ is a homomorphism of $A$--bimodules. We need to
prove that $\gamma$ is a cointegral for $\rcomatrix BM$.
 Given $\varphi \tensor{} m, \varphi' \tensor{} m' \in
\rcomatrix{B}{M}$,
\begin{eqnarray*}
\sum_i (\varphi \tensor{} e_i)\gamma(\rdual{e_i} \tensor{} m
\tensor{} \varphi' \tensor{} m') & = &\sum_i (\varphi \tensor{}
e_i)\rdual{e_i}(\mathsf{s}(m \tensor{} \varphi')(m')) \\
 & = & \varphi \tensor{B} \mathsf{s}(m \tensor{} \varphi')(m') \\
 &= & \varphi \mathsf{s}(m \tensor{} \varphi') \tensor{B} m' \\
 & = & \sum_i\varphi(\mathsf{s}(m \tensor{} \varphi')(e_i))(\rdual{e_i}
 \tensor{} m') \\
 & = & \sum_i\gamma (\varphi \tensor{} m \tensor{} \varphi' \tensor{}
 e_i)(\rdual{e_i} \tensor{} m'),
\end{eqnarray*}
where the second equality follows from the dual basis property.
Furthermore
\begin{eqnarray*}
(\gamma \circ \Delta_{\rcomatrix{B}{M}})(\varphi \tensor{} m) & =
& \sum_i \gamma (\varphi \tensor{} e_i \tensor{} \rdual{e_i}
\tensor{} m)
 =  \sum _i \varphi(\mathsf{s}(e_i \tensor{} \rdual{e_i})(m)) \\
  & = & \varphi(m) = \eps_{\rcomatrix{B}{M}}(\varphi \tensor{} m).
\end{eqnarray*}
Thus $\gamma$ is a cointegral in $\rcomatrix{B}{M}$, i.e., the
comatrix coring $\rcomatrix{B}{M}$ is coseparable as required.

(2) Suppose that  $\rcomatrix{B}{M}$ is coseparable, and let
$\gamma$ be the corresponding cointegral. We  aim to show that
the pre-cointegral $\widetilde{\gamma}$  in the canonical coring
$S\ot_B S$ constructed in Lemma~\ref{cointegral} is a cointegral.
In view of the definition of the coproduct and counit in
Sweedler's coring this is equivalent to showing   that for all $s,
s' \in S$, $\widetilde{\gamma} (s \tensor{} 1_S \tensor{} s') =
ss'$. We freely use the identification of $S$ with
$M\ot_A\rdual{M}$ described in Remark~\ref{rem.s} to compute
\begin{eqnarray*}
\widetilde{\gamma}(s \tensor{} 1_S \tensor{} s') & = &
\sum_{i,j,k}\widetilde{\gamma}(s(e_j) \tensor{} \rdual{e_j}
\tensor{} e_i
\tensor{} \rdual{e_i} \tensor{} s'(e_k) \tensor{} \rdual{e_k}) \\
 & = & \sum_{i,j,k} s(e_j) \gamma(\rdual{e_j} \tensor{} e_i
\tensor{} \rdual{e_i} \tensor{} s'(e_k)) \tensor{} \rdual{e_k} \\
 & = & \sum_{j,k}s(e_j) \gamma (\Delta_{\rcomatrix BM}( \rdual{e_j} \tensor{} s'(e_k)))
 \tensor{} \rdual{e_k} \\
  & = & \sum_{j,k}s(e_j)\rdual{e_j}(s'(e_k)) \tensor{} \rdual{e_k}
   =  \sum_{k}s(s'(e_k)) \tensor{} e_k = ss',
\end{eqnarray*}
as required. Note that the fourth equality follows from the fact
that the pre-cointegral $\gamma$ is a cointegral. Therefore we
conclude that $\widetilde{\gamma}$ is a cointegral for $S\ot_B S$,
i.e., the canonical coring is coseparable as asserted.
\end{proof}

Note that Theorem~\ref{prop.cosep} implies in particular that if
${}_BM_A$ is a separable bimodule, then $S\ot_B S$ is a
coseparable coring. This also follows from
\cite[Theorem~3.1(1)]{Kad:sep} and \cite[26.10]{BrzWis:cor} (the
latter is a refinement of \cite[Corollary~3.7]{Brz:str}).
Theorem~\ref{prop.cosep} leads to a more complete description of
the relationship between separable bimodules and coseparable
comatrix corings in the case of a faithfully flat extension $B\to
S$.

\begin{corollary}\label{cor.cosep}
Let ${}_BM_A$ be a bimodule such that $M_A$ is a finitely
generated projective, and let $S = \rend{A}{M}$. If either ${}_BS$
or $S_B$ is faithfully flat then the following statements are
equivalent
\begin{enumerate}[(a)]
\item $M$ is a separable bimodule.
\item The comatrix coring $\rcomatrix{B}{M}$ is a
coseparable $A$-coring.
\item  $S \tensor{B} S$ is a coseparable $S$-coring.
\end{enumerate}
\end{corollary}
\begin{proof} The implications $(a)\Rightarrow (b) \Rightarrow (c)$ are contained in
Theorem~\ref{prop.cosep}. Suppose that $S \tensor{B} S$ is a coseparable $S$-coring.
In view of the faithful flatness, $B\to S$ is a split extension
by \cite[Corollaries~3.6,~3.7]{Brz:str}. Since $M_A$ is finitely generated projective,
\cite[Theorem~3.1(2)]{Kad:sep} implies that $M$ is a separable bimodule.
This  proves the implication $(c)\Rightarrow (a)$, and completes the proof of the corollary.
\end{proof}

Note that $S$ is a faithfully flat left $B$-module if $M$ is a
faithfully flat left $B$-module.

Recall from \cite{Mor:adj} that a ring extension $B\to S$ is a
Frobenius extension if and only if the restriction of scalars
functor has the same right and left adjoint (cf.\
\cite{Mor:adj}). Following this observation a functor is called a
{\em Frobenius functor} in case it
 has the same right and left adjoint (cf.\  \cite{CaeMil:Doi},
\cite{CasGom:Fro}). Motivated by this correspondence between
Frobenius extensions and Frobenius functors one says that  an
$A$-coring $\cC$ is {\em Frobenius} provided the forgetful functor
from the category of right $\cC$-comodules to the category of
right $A$-modules is Frobenius. Equivalently, $\cC$ is a Frobenius
coring if and only if there exist an invariant $e\in \cC^A =
\{c\in \cC \; |\; \forall a\in A, \; ac=ca\}$ and a pre-integral
$\gamma:\cC\otimes_A\cC\to A$ such that for all $c\in \cC$, $
\gamma(c\otimes_{A}e) = \gamma(e\otimes_{A}c) =
    \eC(c)$.
The pair $(\gamma, e)$ is called a {\em reduced Frobenius system}
for $\cC$ \cite{Brz:tow}.

\begin{theorem}\label{prop.Frob}
Let ${}_BM_A$ be a bimodule such that $M_A$ is a finitely
generated projective, and let $S = \rend{A}{M}$. Then
\begin{enumerate}[(1)]
\item If $M$ is a Frobenius bimodule, then $\rcomatrix{B}{M}$ is a
Frobenius $A$--coring.
\item If $\rcomatrix{B}{M}$ is a Frobenius $A$-coring, then $S
\tensor{B} S$ is a Frobenius $S$--coring.
\end{enumerate}
\end{theorem}
\begin{proof}
$(1)$ Let $\coring{C} = \rcomatrix{B}{M}$, and denote by  $R$ the
opposite ring of $\ldual{\cC}$, i.e. $R =
(\ldual{\coring{C}})^{opp}$. In view of  \cite[Theorem
4.1]{Brz:str}, to prove that $\rcomatrix{B}{M}$ is a Frobenius
coring suffices it to construct an $(A,R)$-bimodule isomorphism
$\coring{C} \cong R$. On the other hand, by \cite[Proposition
1]{KaoGom:com}, there is a ring isomorphism $R \cong \lend{B}{M}$,
such  that the right $R$-module structure on $\coring{C}$ is given
by $(\varphi \tensor{} m) \cdot r = \varphi \tensor{} r(m)$ for
$\varphi \tensor{} m \in \coring{C}$ and $r \in R$ viewed as an
element of $\lend{B}{M}$. Suppose that $M$ is a Frobenius bimodule
and let   $\theta : \rdual{M} \rightarrow \ldual{M}$ be the
(defining) Frobenius $(A,B)$-bimodule isomorphism.  Define an
$(A,B)$-bimodule isomorphism $\iota : \coring{C} \cong R$ as the
composite
\[
\xymatrix{\rcomatrix{B}{M} \ar^-{\theta \tensor{} M}[rr] & &
\ldual{M} \tensor{B} M \cong \lend{B}{M},}
\]
where the last isomorphism follows from the fact that ${}_BM$ is a
finitely generated projective module. The isomorphism $\iota$
explicitly comes out as $\iota (\varphi \tensor{} m)(x) = \theta
(\varphi)(x)m$ for $\varphi \tensor{} m \in \coring{C}$, $x \in
M$. A routine calculation verifies that $\iota$ is an
$(A,R)$--bimodule map.

$(2)$ Suppose that $\rcomatrix{B}{M}$, and let  $\gamma :
\rcomatrix{B}{M} \tensor{A} \rcomatrix{B}{M} \rightarrow A$ and
$e = \sum_a\rdual{w_a} \tensor{} w_a \in (\rcomatrix{B}{M})^A$ be
a reduced Frobenius system. This means that $\gamma$ is a
pre-cointegral and  for all $\varphi \tensor{} m \in
\rcomatrix{B}{M}$
$$
\sum _a \gamma (\varphi \tensor{} m \tensor{} \rdual{w_a}
\tensor{} w_a) = \sum_a \gamma (\rdual{w_a} \tensor{} w_a
\tensor{} \varphi \tensor{} m) = \varphi (m). \eqno{(*)}
$$
Consider the pre-cointegral  $\widetilde{\gamma} : S \tensor{B} S
\tensor{B} S \rightarrow S$ constructed in Lemma \ref{cointegral},
and define $\widetilde{e} = \sum_{i,a}e_i \tensor{} \rdual{w_a}
\tensor{} w_a \tensor{} \rdual{e_i} \in (S \tensor{B} S)^S$. In
this definition and throughout the rest of the proof we freely use
the identification of $S$ with $M\ot_A\rdual{M}$ described in
Remark~\ref{rem.s}. We need to check that
$(\widetilde{\gamma},\widetilde{e})$ is a Frobenius system for
$S\ot_B S$, i.e.,   that
\[
\widetilde{\gamma}(s \tensor{} s' \widetilde{e}) = ss' =
\widetilde{\gamma}(\widetilde{e} s \tensor{} s') ,
\]
for all $s, s' \in S$. This is carried out by  the following
explicit computations. First,
\begin{eqnarray*}
\widetilde{\gamma}(s \tensor{} s' \widetilde{e}) &=&
\sum_{i,a}\widetilde{\gamma}(s \tensor{} s'(e_i) \tensor{}
\rdual{w_a} \tensor{} w_a \tensor{} \rdual{e_i})
\\
& = &
 \sum_{i,k,a}\widetilde{\gamma}(s(e_k) \tensor{} \rdual{e_k} \tensor{}
 s'(e_i) \rdual{w_a} \tensor{} w_a \tensor{} \rdual{e_i}) \\
 & = & \sum_{i,k,a} s(e_k)\gamma(\rdual{e_k} \tensor{} s'(e_i) \tensor{}
 \rdual{w_a} \tensor{} w_a) \tensor{} \rdual{e_i} \\
  &= & \sum_{i,k}s(e_k)\rdual{e_k}(s'(e_i)) \tensor{} \rdual{e_i}
  =  \sum_{i,k}(s(e_k) \tensor{} \rdual{e_k})(s'(e_i) \tensor{}
  \rdual{e_i})
  = ss' ,
\end{eqnarray*}
as required. Note that the first equality follows from the
definition of $\widetilde{e}$ and already incorporates the formula
for the product in $S$ in terms of elements of $M\ot_A\rdual{M}$
as explained in Remark~\ref{rem.s}. Remark~\ref{rem.s} is also
used to derive the second and fifth equalities. The penultimate
equality follows from equation ($*$).  Second
\begin{eqnarray*}
ss'   & = & \sum_{k}s(e_k) \tensor{} \rdual{e_k}s'
  = \sum_{i,k} e_i \rdual{e_i}s(e_k) \tensor{} \rdual{e_k}s' \\
  & = & \sum_{i,k}e_i \tensor{} \rdual{e_i}s(e_k)\rdual{e_k}s'
  =  \sum_{i,k,a}e_i \tensor{} \gamma (\rdual{w_a} \tensor{} w_a \tensor{}
  \rdual{e_i}s \tensor{} e_k)\rdual{e_k}s' \\
  & = & \sum_{i,k,a}\widetilde{\gamma}(e_i \tensor{} \rdual{w_a} \tensor{} w_a
  \tensor{} \rdual{e_i}s \tensor{} e_k \tensor{} \rdual{e_k}s')
= \widetilde{\gamma}(\widetilde{e} s \tensor{} s'),
\end{eqnarray*}
as required. Here Remark~\ref{rem.s} is used in derivation of the
first two and the last equalities, while the fourth equality
follows from equation ($*$). Thus we have proven that
$(\widetilde{\gamma},\widetilde{e})$ is a Frobenius system for
$S\ot_B S$, so that $S\ot_B S$ is a Frobenius coring as asserted.
\end{proof}

Note that Theorem~\ref{prop.Frob} implies in particular that if
${}_BM_A$ is a Frobenius bimodule, then $S\ot_B S$ is a Frobenius
coring. This also follows from the endomorphism ring theorem
\cite[Theorem~2.5]{Kad:new} and \cite[Theorem~2.7]{Brz:tow}. As
was the case for coseparable comatrix corings,
Theorem~\ref{prop.Frob} leads to a more complete description of
the relationship between Frobenius bimodules and Frobenius
comatrix corings in the case of a faithfully flat extension $B\to
S$ which in addition satisfies a weak version of Williard's
condition.

\begin{corollary}\label{cor.Frob}
Let ${}_BM_A$ be a bimodule such that both $M_A$ and ${}_BM$ are
finitely generated projective, and let $S = \rend{A}{M}$. Suppose
that either ${}_BS$ or $S_B$ is faithfully flat  and that $\lhom
SMS \cong \rhom AMA$ as $(A,B)$-bimodules. Then the following
statements are equivalent
\begin{enumerate}[(a)]
\item $M$ is a Frobenius bimodule.
\item The comatrix coring $\rcomatrix{B}{M}$ is a
Frobenius $A$-coring.
\item  $S \tensor{B} S$ is a Frobenius $S$-coring.
\end{enumerate}
\end{corollary}
\begin{proof} The implications $(a)\Rightarrow (b) \Rightarrow (c)$ are contained in
Theorem~\ref{prop.Frob}. Suppose that $S \tensor{B} S$ is a Frobenius $S$-coring. In
view of the faithful flatness, $B\to S$ is a Frobenius extension by \cite[Theorem~2.7]{Brz:tow}.
Since $\lhom SMS \cong \rhom AMA$ as $(A,B)$-bimodules, the converse of the endomorphism ring theorem
\cite[Theorem~2.8]{Kad:new} implies that $M$ is a Frobenius bimodule. This  proves the implication
$(c)\Rightarrow (a)$, and completes the proof of the corollary.
\end{proof}

As noted in \cite[Section~2.3]{Kad:new}, the condition $\lhom SMS
\cong \rhom AMA$ as $(A,B)$-bimodules is in particular satisfied
when $M_A$ is a generator module.

\begin{remark}
The central idea of this paper is that  properties of a bimodule ${}_BM_A$ imply
analogous properties of the endomorphism ring $S = \rend{A}{M}$. These in turn lead to 
 corresponding properties of   the  Sweedler
$S$--coring associated to the extension $B\to S$. The comatrix $A$--coring built with $M$ can be thought of as a dual of the endomorphism ring, and thus can be envisioned as lying  in between a bimodule $M$ and the Sweedler coring associated to
$B\to S$. Thus, combining the results of the present paper with that of existing literature, the situation can be summarised in terms of  the following deductive diagrams.

 In case $M_A$ is finitely generated and projective, 
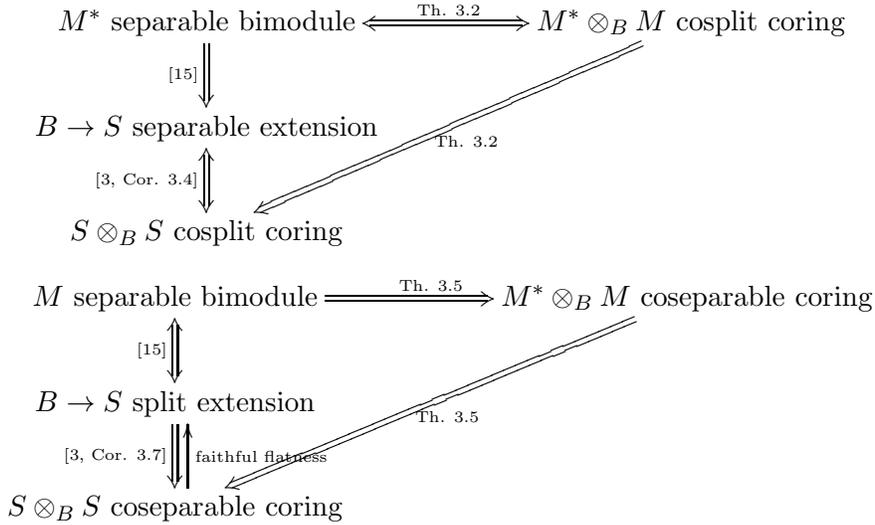
\[
\xymatrix{ M^*  \hbox{ separable bimodule}\ar@{<=>}^{\hbox{\tiny Th. \ref{th.cosplit}}}[rr]
\ar@{=>}_{\hbox{\tiny \cite{Sug:not}}}[d]& & M^* \tensor{B} M \hbox{ cosplit coring} \ar@{=>}^{\hbox{\tiny Th. \ref{th.cosplit}}}[ddll]\\
B \rightarrow S \hbox{ separable extension} \ar@{<=>}_{\hbox{\tiny \cite[Cor. 3.4]{Brz:str}}}[d]& &  \\
S \tensor{B} S \hbox{ cosplit coring} & &    }
\]
\[
\xymatrix{M  \hbox{ separable bimodule}\ar@{=>}^{\hbox{\tiny Th. \ref{prop.cosep}}}[rr]
\ar@{<=>}_{\hbox{\tiny \cite{Sug:not}}}[d]& & M^* \tensor{B} M \hbox{ coseparable coring} \ar@{=>}^{\hbox{\tiny Th. \ref{prop.cosep}}}[ddll]\\
B \rightarrow S \hbox{ split extension} \ar@{=>}_{\hbox{\tiny \cite[Cor. 3.7]{Brz:str}}}[d] & &  \\
S \tensor{B} S \hbox{ coseparable coring}  \ar@<-1ex>_{\hbox{\tiny
faithful flatness }}[u] & & }
\]

In case ${}_BM$ and $M_A$ are finitely generated and projective, 
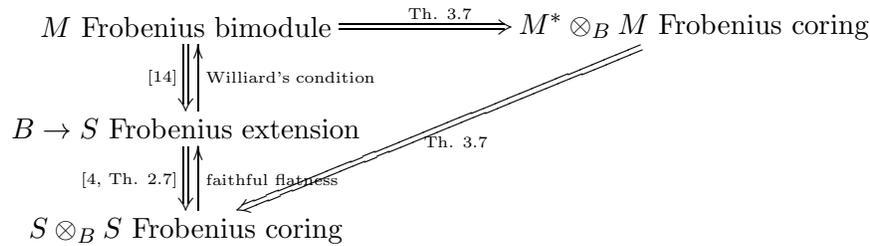
\[
\xymatrix{M  \hbox{ Frobenius bimodule}\ar@{=>}^{\hbox{\tiny Th. \ref{prop.Frob}}}[rr]
\ar@{=>}_{\hbox{\tiny \cite{Mor:adj}}}[d]& & M^* \tensor{B} M \hbox{ Frobenius coring} \ar@{=>}^{\hbox{\tiny Th. \ref{prop.Frob}}}[ddll]\\
B \rightarrow S \hbox{ Frobenius extension} \ar@{=>}_{\hbox{\tiny
\cite[Th. 2.7]{Brz:tow}}}[d] \ar@<-1ex>_{\hbox{\tiny
Williard's condition }}[u]& &  \\
S \tensor{B} S \hbox{ Frobenius coring} \ar@<-1ex>_{\hbox{\tiny
faithful flatness }}[u] & &    }
\]
\end{remark}

\section*{Acknowledgements}
Tomasz Brzezi\'nski  would like to thank the Engineering and
Physical Sciences Research Council for an Advanced Fellowship. He
also thanks the Department of Algebra, University of Granada for
hospitality. \\
Investigaci\'{o}n parcialmente financiada por el Proyecto BFM2001-3141
del Ministerio de Ciencia y Tecnolog\'{\i}a de Espa\~{n}a.


\begin{thebibliography}{99}

\parindent0pt
\footnotesize

\bibitem{AndFul:rin} F.\ Anderson and K.\ Fuller. {\em Rings and Categories of
    Modules}, Springer, Berlin, 1974.

\bibitem{Bas:KTh} H.\ Bass. {\em Algebraic K-Theory}. W.A.\
    Benjamin, inc., New York, 1968.
\bibitem{Brz:str} T.\ Brzezi\'nski. The structure of corings. Induction functors,
Maschke-type theorem, and Frobenius and
        Galois-type properties. {\em Alg.\
    Rep.\ Theory}, 5: 389--410, 2002.
\bibitem{Brz:tow} T.\ Brzezi\'nski. Towers of corings.  {\em Preprint} ArXiv: math.RA/0201014, 2002. {\em Commun.\ Algebra} to appear.
\bibitem{BrzKad:bis} T.\ Brzezi\'nski,  L.\ Kadison and R.\ Wisbauer. On
      coseparable and biseparable  corings.
      {\em Preprint} ArXiv math.RA/0208122, 2002.
\bibitem{BrzWis:cor} T.\ Brzezi\'nski and R.\ Wisbauer. {\em Corings and Comodules}, Cambridge University Press, Cambridge, to appear.
\bibitem{CaeKad:bis} S.\ Caenepeel and L.\ Kadison. Are biseparable extensions Frobenius?
{\em   K-Theory}, 24: 361--383, 2001.
\bibitem{CaeMil:Doi}  S.\ Caenepeel, G.\ Militaru and  S.\ Zhu.
    Doi-Hopf modules, Yetter-Drinfel'd modules and Frobenius type
    properties.  {\em Trans.\ Amer.\ Math.\ Soc.}  349:4311--4342, 1997.
    \bibitem{CasGom:Fro} F.\ Casta\~no Iglesias, J.\  G\'omez-Torrecillas and
     C.\ N\u ast\u asescu. Frobenius functors.
        Applications.  {\em Commun.\ Algebra},  27:4879--4900, 1999.
\bibitem{KaoGom:com} L.\ El Kaoutit and  J.\ G\'omez-Torrecillas.
 Comatrix corings: Galois corings, descent theory, and a structure
theorem for cosemisimple corings, {\em Preprint} ArXiv:
math.RA/0207205, 2002. To appear in {\em Math.\ Z.}
\bibitem{Guz:coi} F. Guzman. Cointegrations, relative cohomology for
  comodules and coseparable corings. {\em J.\ Algebra}, 126:211--224, 1989.
\bibitem{Kad:new}  L.\ Kadison. {\em New Examples of Frobenius
    Extensions}, AMS, Providence R.I., 1999.
\bibitem{Kad:sep}  L.\ Kadison. Separability and the twisted Frobenius bimodule.
{\em   Alg.\
    Rep.\ Theory}, 2: 397--414, 1999.
\bibitem{Mor:adj}  K.\ Morita. Adjoint pairs of functors and
        Frobenius extensions. {\em Sci.\ Rep.\ Tokyo Kyoiku Daigaku
        Sect.\ A}, 9:40--71, 1965.
\bibitem{Sug:not} K.\ Sugano. Note on separability of endomorphism rings. {\em Hokkaido Math.\ J.}, 11:111--115, 1982.
\bibitem{Swe:pre}
M.\ Sweedler. The predual theorem to the Jacobson-Bourbaki
theorem, {\em Trans.\ Amer.\ Math.\ Soc.} 213:391--406, 1975.
\end{thebibliography}
\end{document}